\documentclass{amsart}
\usepackage{amsmath,amssymb,bm}
\usepackage{latexsym}

\title
{Reeb orbits trapped by Denjoy minimal sets}
\author{Takahiro Arai, Takashi Inaba and Yosuke Kano}
\address{
Department of Mathematics and Informatics,
Graduate School of Science,
Chiba University,
Chiba 263-8522,
Japan
}
\email{inaba@math.s.chiba-u.ac.jp}

\newtheorem{Thm}{Theorem}[section]
\newtheorem{Lem}[Thm]{Lemma}
\newtheorem{Prop}[Thm]{Proposition}
\newtheorem{Cor}[Thm]{Corollary}

\newcommand{\R}{\mathbb{R}}

\bmdefine{\boldr}{r}
\bmdefine{\boldv}{v}

\makeatletter
    
    \@addtoreset{equation}{section}
  \makeatother

\begin{document}
\maketitle
\begin{abstract}
Let $\varphi$ be any flow on $T^n$
obtained as the suspension of a diffeomorphism of $T^{n-1}$
and let $\mathcal A$ be any compact invariant set of $\varphi$.
We realize $(\mathcal A, \varphi|_{\mathcal A})$ 
up to
reparametrization
as an invariant set of the Reeb flow 
of a contact form on $\R^{2n+1}$ 
equal to the standard contact form outside a compact set
and defining the standard contact structure on all of  $\R^{2n+1}$.
This generalizes the construction of Geiges, R\"ottgen and Zehmisch.
\end{abstract}

\section{Introduction}
In \cite{EH} Eliashberg and Hofer proved that
if $\alpha$ is a contact form on $\R^3$ which 
coincides with the standard contact form outside a compact set
and defines the standard contact structure on all of  $\R^3$
and if the Reeb flow of $\alpha$ has a bounded forward orbit,
then the flow necessarily has a periodic orbit. 
Recently, Geiges, R\"ottgen and Zehmisch \cite{GRZ} (see also R\"ottgen \cite{R})
showed that the higher dimensional analogue 
of the Eliashberg-Hofer theorem is not true.
They constructed, for any $n\ge2$, a contact form on $\R^{2n+1}$ equal
to the standard contact form outside a compact set
and defining the standard contact structure on whole $\R^{2n+1}$,
with the property that its Reeb flow has a bounded forward orbit 
but has no periodic orbits.
In their example bounded forward orbits are trapped by 
an invariant $n$-dimensional torus with an ergodic linear flow.

In this note we generalize their construction.
Our result enables us to produce diverse compact invariant sets
in Reeb flows. 
Let  $(\theta_1,\cdots,\theta_n)$ be the standard angle coordinates of 
the $n$-dimensional torus $T^n$.
We denote the coordinates of $\R^{2n}$ 
(resp. $\R^{2n+1}$)
by $(x_1,y_1,\cdots,x_n,y_n)$ (resp. $(x_1,y_1,\cdots,x_n,y_n,z)$),
and the polar coordinates in the $(x_j,y_j)$-plane by $(r_j,\theta_j)$.
Via a natural inclusion using polar coordinates,
we identify $T^n$ with the subset $\mathbb T=\{r_1=\cdots =r_n=1\}$
of $\R^{2n}$.

We show the following

\begin{Thm}\label{main}
Let $V$ be any vector field on $\mathbb T$
satisfying  $\sum_{j=1}^nd\theta_j(V)>0$,
and let $\mathcal A\subset \mathbb T$ be any compact invariant set of
the flow generated by $V$.
Then, one can find a $C^\infty$ contact form $\alpha$ on $\R^{2n+1}$
which defines the standard contact structure
$\xi_{\mathrm{st}}$ and satisfies the following properties:
\begin{itemize}
\item[$(1)$] $\alpha$ equals the standard contact form $\alpha_{\mathrm{st}}$
outside a compact neighborhood of $\mathcal A\times\{0\}$.
\item[$(2)$] The flow $\varphi$ generated by the Reeb vector field $R$ 
of $\alpha$ has $\mathcal A\times\{0\}$ as its invariant set,
and there exists a positive $C^\infty$ function $f$ on $\mathbb T$
such that $R=f V$ on $\mathcal A\times\{0\}$.
\item[$(3)$] All orbits of $\varphi$ in the complement of 
$\mathcal A\times\{0\}$ in $\R^{2n+1}$
are unbounded.
\item[$(4)$] There exists an orbit of $\varphi$ 
which is bounded in forward time
and unbounded in backward time.
\end{itemize}
\end{Thm}

\noindent\textbf{Remark.}
On $\mathbb T$,  the plane field $\xi_{\mathrm{st}}\cap T(\mathbb T)$ is
expressed by $\sum_{j=1}^nd\theta_j=0$
 (See \S 2 for the definition of $\xi_{\mathrm{st}}$).
Thus, $V$ in the above theorem covers
all vector fields  on  $\mathbb T$ positively transverse to $\xi_{\mathrm{st}}$.
\vspace{3mm}

Let $\psi$ be a flow on a manifold $M$.
We say in this note that a compact minimal set $\mathcal M$ of $\psi$
is a \textit{Denjoy minimal set}
if for any point $p$ of $\mathcal M$
there exists a codimension $1$ disk $D$ in $M$ transverse to $\psi$ 
passing through $p$
such that $D\cap \mathcal M$ is homeomorphic to a Cantor set.
Combining Hall's result \cite{H} with our theorem,
we obtain the following

\begin{Cor}\label{Denjoy}
For any $n\ge3$,
there exists a $C^\infty$ contact form $\alpha$ on
$\R^{2n+1}$ which defines the standard contact structure $\xi_{\mathrm{st}}$ 
and satisfies the following properties:
\begin{itemize}
\item[$(1)$]  $\alpha$ equals $\alpha_{\mathrm{st}}$
outside a compact set.
\item[$(2)$] The flow $\varphi$ generated by the Reeb vector field
of $\alpha$ admits a Denjoy minimal set $\mathcal M$.
And $\mathcal M$ is the unique compact invariant set for $\varphi$
$($in particular, there are no periodic orbits$)$.
\item[$(3)$] There exists an orbit of $\varphi$
which is bounded in forward time
and unbounded in backward time
$($in this case the $\omega$-limit set of this orbit is $\mathcal M$
$)$.
\end{itemize}
\end{Cor}

\section{Basic definitions and facts}
A $C^\infty$ $1$-form $\alpha$ on a $(2n+1)$-dimensional manifold $M$
is a {\it contact form} if 
$\alpha\wedge(d\alpha)^n$ never vanishes on $M$.
A codimension $1$ plane field $\xi$ on $M$ is a 
{\it contact structure} on $M$
if $\xi=\operatorname{ker}\alpha$ for
some contact form $\alpha$.
A manifold with a contact structure is called a {\it contact manifold}.
Given a contact form $\alpha$ the {\it Reeb vector field} of $\alpha$ is
the unique vector field $R$ satisfying $\alpha(R)=1$ and $i(R)d\alpha=0$.
For a contact manifold $(M,\xi)$
a vector field $X$ on $M$ is a {\it contact vector field} of $\xi$
if the flow generated by $X$ preserves $\xi$.
We denote by $C^\infty(M)$ the set of all real valued $C^\infty$
functions on $M$, by $\Gamma^\infty(\xi)$
the set of all $C^\infty$ vector fields on $M$
tangent to $\xi$ at each point of $M$,
and by
$\Gamma_{\xi}^\infty(TM)$
the set of all
$C^\infty$ contact vector fields on $M$.

The following two facts are fundamental (see \cite{G}).
\begin{Prop}\label{corres}
Let $\xi$ be a $C^\infty$ contact structure on $M$
and $\alpha$ a $C^\infty$ contact form such that
$\xi=\operatorname{ker}\alpha$.
Then,
there is a bijective correspondence
between $\Gamma_{\xi}^\infty(TM)$ and $C^\infty(M)$
given as follows:
Define $\Phi:\Gamma_{\xi}^\infty(TM)\to C^\infty(M)$
by $\Phi(X)=\alpha(X)$
and $\Psi:C^\infty(M)\to\Gamma_{\xi}^\infty(TM)$ by
$\Psi(H)=HR+Y$,
where $R$ is the Reeb vector field of $\alpha$
and $Y$is the unique vector field in $\Gamma^\infty(\xi)$
such that $i(Y)d\alpha=dH(R)\alpha-dH$.
Then $\Phi$ and $\Psi$ are the inverses of each other.
\end{Prop}

\begin{Prop}\label{reeb}
Let $\alpha$ be a $C^\infty$ contact form on $M$
and let $\xi=\operatorname{ker}\alpha$.
If $X\in\Gamma_{\xi}^\infty(TM)$ and
if $X$ is transverse to $\xi$ everywhere on $M$,
then $X$ is the Reeb vector field of the 
contact form $\alpha/\alpha(X)$.
\end{Prop}

We call
\begin{equation}
\alpha_{\mathrm{st}}
=dz+\frac12\sum_{j=1}^{n}
(x_jdy_j-y_jdx_j)
 =dz+\frac12\sum_{j=1}^{n}
{r_j}^2d\theta_j
\end{equation}
the {\it standard contact form} on
$\R^{2n+1}$
and 
$\xi_{\mathrm{st}}=\operatorname{ker}\alpha_{\mathrm{st}}$
the {\it standard contact structure} on 
$\R^{2n+1}$.
The co-orientation of $\xi_{\mathrm{st}}$ is given by $\alpha_{\mathrm{st}}$.
It is easy to see that the Reeb vector field of 
$\alpha_{\mathrm{st}}$ is $\partial_z$.
The vector field $Y$ in Proposition \ref{corres}
with respect to $\alpha_{\mathrm{st}}$ is
expressed as follows
(at points with $r_j\ne0$).
\begin{equation}
Y=\sum\left[\left(\frac{r_j}2H_z-\frac{H_{\theta_j}}{r_j}\right)
\partial_{r_j}
+\frac{H_{r_j}}{r_j}\left(\partial_{\theta_j}-
\frac{{r_j}^2}2\partial_z\right)\right]
\label{Y}
\end{equation}
Here,
$\displaystyle
\left\{
\partial_{r_j}, 
\partial_{\theta_j}-
\dfrac{{r_j}^2}2\partial_z
\right\}_j$
is a frame of $\xi_{\mathrm{st}}$.

\section{Conditions on $X$ and $H$}

From now on, we exclusively consider
the standard contact structure $\xi_{\mathrm{st}}$ on $\R^{2n+1}$.
Take 
$C^\infty$ functions $k_j$ ($j=1,\cdots, n$) on $T^n$
satisfying the condition
\begin{equation}\label{k}
\sum_{j=1}^nk_j=1.
\end{equation}
We denote by $U$ the subset of 
$\R^{2n}$ consisting of points
with $r_1r_2\cdots r_n\ne0$.
Then,
by composing $k_j$ with
the natural projection $U\to T^n$ in polar coordinates,
we may regard $k_j$ as functions on $U$
(although in general we cannot extend the domain to whole $\R^{2n}$).
Let $\mathcal A\subset T^n$ be a compact invariant set for the flow
generated by the vector field
$\sum_{j=1}^nk_j\partial_{\theta_j}$.
Let $\lambda>0$ be a constant.

We consider the following conditions for 
a $C^\infty$ contact vector field $X$ of 
$\xi_{\mathrm{st}}$.
\begin{itemize}
\item[(X1)] $X$ is positively transverse to $\xi_{\mathrm{st}}$
at each point of $\R^{2n+1}$.
\item[(X2)] $X=\partial_z$ outside a compact neighborhood of 
$\mathcal A\times\{0\}$.
\item[(X3)] $dz(X)>0$ outside $\mathcal A\times\{0\}$.
\item[(X4)] $X=2\lambda\sum_{j=1}^nk_j\partial_{\theta_j}$
on $\mathcal A\times\{0\}$.
\end{itemize}
Next, we consider the following conditions for 
a $C^\infty$ function $H$ on $\R^{2n+1}$.
\begin{itemize}
\item[(H1)] $H$ is positive at each point of $\R^{2n+1}$.
\item[(H2)] $H=1$ outside a compact neighborhood of 
$\mathcal A\times\{0\}$.
\item[(H3)] $\displaystyle
H-\dfrac{\ 1\ }2\sum_jr_jH_{r_j}>0$ outside $\mathcal A\times\{0\}$.
\item[(H4)] On $\mathcal A\times\{0\}$, we have
$H=\lambda$, 
$H_{r_j}=2\lambda k_j$, $H_{\theta_j}=0$ ($j=1,\cdots, n$), 
and
$H_z=0$.
\end{itemize}
Then, via the correspondence (Proposition \ref{corres})
these two sets of conditions are closely connected with each other:

\begin{Lem}\label{X and H}
Suppose a $C^\infty$ function $H$ on $\R^{2n+1}$
satisfies the conditions $(${\rm H1}$)$ to $(${\rm H4}$)$.
Then,
the contact vector field $X$ of $\xi_{\mathrm{st}}$
which corresponds with $H$ via $\alpha_{\mathrm{st}}$
satisfies the conditions $(${\rm X1}$)$ to $(${\rm X4}$)$.
\end{Lem}

\noindent\textbf{Proof.}
By Proposition \ref{corres},
(H1) implies (X1).
Suppose $H$ satisfies (H2) with a compact neighborhood $W$.
Then,
since all the partial derivatives of $H$ vanish outside $W$,
by (\ref{Y}) we have $Y=0$ hence
$X=H\partial_z+Y=\partial_z$ outside $W$.
This shows (X2).
By (\ref{Y}) we have
$\displaystyle
dz(X)=dz(H\partial_z+Y)=H+dz(Y)=
H-\dfrac{\ 1\ }2\sum r_jH_{r_j}$.
Thus, (X3) implies (X3).
Finally, suppose $H$ satisfies (H4).
Then (X4) follows from (\ref{Y}) and (\ref{k})
by a simple computation. 
\hfill $\Box$
\vspace{5mm}

\section{Constructing $H$}
This section is devoted to the construction of a function $H$ on
$\R^{2n+1}$ satisfying the conditions (H1) to (H4).
Let $k_j$, $\mathcal A$ and $U$ be as in the preceding section.
For a real number $b$ and $\tau$
we consider the following function $Q^b[\tau]$
on $U$.
\begin{equation}
Q^b[\tau](r_1,\cdots,r_n,\theta_1,\cdots,\theta_n)=
\sum_ik_i(\theta_1,\cdots,\theta_n)(r_i-\tau)^2+b\sum_{p\ne q}(r_p-r_q)^2.
\end{equation}
Let $\boldr$ be the column vector with entries $r_1,\cdots,r_n$
and $A(b)$ the $n\times n$ matrix with
$(i,i)$-th entry $k_i+2(n-1)b$ 
and $(i,j)$-th entry $-2b$ for $i\ne j$.
Then $Q^b[0]={}^t\boldr A(b)\boldr$.

\begin{Lem}\label{posdef}
If $b$ is sufficiently large,
then $Q^b[\tau]$ is non-negative everywhere on $U$.
\end{Lem}

\noindent\textbf{Proof.}
It suffices to show that the matrix $A(b)$ is positive definite.
If $b$ is sufficiently large,
the $s\times s$ principal minor ($s\leq n-1$)
of $A(b)$ is positive
because it is a polynomial in $b$ of degree $s$
with the coefficient of the highest order term being a positive constant.
Also, $\det A(b)$  is positive
because by (\ref{k}) it is a polynomial in $b$ of degree $n-1$
with the coefficient of the highest order term being a positive constant.
\hfill $\Box$
\vspace{5mm}

From now on, we always take $b$ large
so that the conclusion of Lemma \ref{posdef} holds. 
Here, we briefly observe some properties of $A(b)$ from linear algebra.
Let us denote the eigenvalues of $A(b)$ by $\lambda_1(b),\cdots,\lambda_n(b)$.
These are functions on $T^n$.
Since 
$\dfrac{k_i}{2b}$ tends to $0$ uniformly on $T^n$ as $b\to\infty$,
the matrix $\dfrac1{2b}A(b)$ converges uniformly  to the matrix
with diagonal entries $n-1$ and the other entries $-1$,
whose eigenvalues are $0$ and $n$.
The eigenspace $W_0$ corresponding to $0$ is a $1$-dimensional subspace of $\R^n$
generated by the vector $(1,\cdots,1)$
and the eigenspace $W_n$ corresponding to  $n$ is the orthogonal complement
of $W_0$.
Thus, renumbering the eigenvalues
 if necessary,
we may assume that as $b\to\infty$
the ratios of the eigenvalues of $A(b)$ converge uniformly as follows:
$\lambda_1(b)/\lambda_i(b)\to0$ ($i>1$)C
$\lambda_i(b)/\lambda_j(b)\to1$ ($i, j>1$),
and that, in the Grassmanian space, 
the eigenspace corresponding to $\lambda_1(b)$ converges
to $W_0$ and
the sum of the eigenspaces corresponding to $\lambda_i(b)$ ($i\geq2$)
converges to $W_n$.

For any $c>0$,
let $E(c)$ be the subset 
$\{\boldv\in\R^n \mid {}^t\boldv A(b)\boldv \leq c\}$ of $\R^n$
and
$J(c)$
the closed line segment
joining two points
$\pm(\sqrt c,\cdots,\sqrt c)$.
We note that $E(c)$ depends on $(\theta_1,\cdots,\theta_n)\in T^n$.
Since we are assuming that $A(b)$ is positive definite,
$E(c)$ is compact.

\begin{Lem}\label{converge}
As $b$ goes to infinity, 
$E(c)$ converges to $J(c)$
uniformly on $T^n$
with respect to the Hausdorff distance.
\end{Lem}

\noindent\textbf{Proof.}
Let $\ell$ be the $1$-dimensional linear subspace of $\R^n$
generated by $(1,\cdots,1)$.
Then, by (\ref{k}) we have $E(c)\cap\ell=J(c)$
for any $b$.
This with the observation made above implies that 
$E(c)$ converges uniformly
to $J(c)$. 
\hfill $\Box$
\vspace{5mm}

We write $L$ for the subset of $U$
consisting of all points satisfying the following condition:
there exists $i$ ($1\leq i\leq n$) such that
$\frac13<r_i<\frac23$ and that $\frac13<r_j$ for all $j$ ($j\ne i$). 

\begin{Lem}\label{L}
Suppose $0<C<7/9$.
If $b$ is sufficiently large,
then we have $Q^b[2]>1+C$ on $L$. 
\end{Lem}

\noindent\textbf{Proof.}
Let $P^b$ be the domain in $U$
satisfying the condition $Q^b[2]\leq1+C$,
and $J$ 
the subset of $U$
defined by
$2-\sqrt{1+C}\leq r_1=\cdots=r_n\leq2+\sqrt{1+C}$.
Then, by Lemma \ref{converge}, 
$P^b$ converges to $J$ as $b\to\infty$.
Since $\frac23<2-\sqrt{1+C}$,
the closure of $L$ is disjoint from $J$,
hence also from $P^b$ for very large $b$.
This shows the conclusion.
\hfill $\Box$
\vspace{5mm}

Here we prepare two more auxiliary functions.
The first one is a monotone increasing $C^\infty$ function
$\rho:[0,\infty)\to[0,1]$
such that
\begin{itemize}
\item[($\rho$1)] $\rho(r)=0$ if and only if $r\le1/3$, and
\item[($\rho$2)] $\rho(r)=1$  if and only if $r\ge2/3$.
\end{itemize}
The second one is a $C^\infty$ function $\mu$ on $T^n$
such that
\begin{itemize}
\item[($\mu$1)] $\mu=\mu_{\theta_j}=0$ ($j=1,\cdots,n$)
on $\mathcal A$, and
\item[($\mu$2)] $\mu>0$ on $T^n-\mathcal A$.
\end{itemize}
Such a $\mu$ exists.

Now,
let $C$ and $\lambda$ be constants satisfying
 $0<C<7/9$ and $1<\lambda<e^C$,
 and $b$ a constant satysfying the conclusions of 
 Lemmas \ref{posdef} and \ref{L}.
We define a $C^\infty$ fuction $K$ on $\R^{2n+1}$ as follows:
\begin{equation}\label{K}
K=
\lambda
\exp\left\{
\left[
1+C-Q^b[2]-(z^2+\mu)\sum_jr_j
\right]
\prod_{\ell=1}^n\rho(r_\ell)
-C
\right\}
.
\end{equation}
Notice that, although the polar coordinate
functions $r_j$ and $\theta_j$ are
defined only on 
$U\times\R$,
by ($\rho$1) $K$ is well-defined
on whole $\R^{2n+1}$.

We have
\begin{Lem} \label{KKK}
$K$ satisfies the conditions 
$(${\rm H1}$)$, $(${\rm H3}$)$ and $(${\rm H4}$)$.
\end{Lem}

\noindent\textbf{Proof.}
Obviously $K$ is positive everywhere, thus satisfies $(${\rm H1}$)$.
On $\mathcal A\times\{0\}$,
by
(\ref{k}), ($\mu$1) and ($\rho$2),
we immediately obtain that
$K=\lambda$, $K_{\theta_j}=0$ ($j=1,\cdots,n$) and $K_z=0$.
Also, by ($\rho$2), 
we have on $\mathcal A\times\{0\}$,
$
K_{r_j}=
-KQ^b[2]_{r_j}
=2\lambda k_j,
$
showing that $K$ satisfies $(${\rm H4}$)$.
Finally, let us check the condition $(${\rm H3}$)$.
Since $K-(1/2)\sum_jr_jK_{r_j}
=(K/2)(2-\sum_j r_j(\log K)_{r_j})$,
it is enough to show that
$\sum_j r_j(\log K)_{r_j}<2$
outside $\mathcal A\times\{0\}$.
We have
\begin{equation}\label{guy}
\begin{split}
&\sum_j r_j(\log K)_{r_j} \\
&\qquad = 
\left[-\sum_j r_jQ^b[2]_{r_j}-(z^2+\mu)\sum_jr_j\right]
\prod_\ell\rho(r_\ell) \\
&\qquad\quad
+\left[
1+C-Q^b[2]-(z^2+\mu)\sum_jr_j
\right]
\sum_j r_j\rho'(r_j)\prod_{\ell\ne j}\rho(r_\ell)
.
\end{split}
\end{equation}
Here we note that
\begin{equation}\label{dQ}
\sum_j r_jQ^b[2]_{r_j}=2Q^b[1]-2.
\end{equation}
We consider three cases separately.

\noindent
\textit{Case 1. $r_j\le1/3$ for some $j$.}
Then, since $\rho(r_j)=\rho'(r_j)=0$,
the RHS of (\ref{guy}) is $0$.

\noindent
\textit{Case 2. $r_j\ge2/3$ for all $j$.}
Then, the RHS of (\ref{guy})
is $-\sum_j r_jQ^b[2]_{r_j}-(z^2+\mu)\sum_jr_j$.
This with ($\mu 2$), (\ref{dQ}) and the positive definiteness of $Q$,
implies the desired conclusion.

\noindent
\textit{Case 3. $1/3<r_j<2/3$ for some $j$.}
In this case, again by  (\ref{dQ}) and the positive definiteness of $Q$
we see that the first term of the RHS of (\ref{guy}) is less than $2$.
On the other hand, Lemma \ref{L} implies that the second term
is non-positive. 
So the conclusion follows.

\noindent
This finishes the proof that $K$ satisfies (H3).
\hfill $\Box$
\vspace{5mm}

\begin{Lem}\label{KeC}
$K$ is not greater than $\lambda e^{-C}$
outside some compact region.
\end{Lem}

\noindent\textbf{Proof.}
In the case where $r_j\le1/3$ for some $j$,
we have $K=\lambda e^{-C}$.
Let us consider the case where 
$r_j>1/3$ for all $j$.
We denote by $T$ the quantity in the square bracket
in (\ref{K}) and claim that
$T$ is non-positive outside some compact region.
In fact, since $Q^b[2]$ and $\mu$ are non-negative and $\sum_jr_j>n/3$,
$T$ is non-positive for $z^2\geq\dfrac{3(1+C)}{n}$.
$T$ is also non-positive for $||\boldr||$ large
because the proof of Lemma \ref{L}
implies that
$\sup\{||\boldr|| \mid Q^b[2](\boldr)\leq1+C\}$ is finite.
This shows the claim,
and the Lemma follows.
\hfill $\Box$
\vspace{5mm}

$K$ constructed above does not satisfy (H2).
By a modification we will improve $K$ 
so as to satisfy all the required conditions.
To this end we need one more auxiliary fuction.
It is a monotone increasing
$C^\infty$ function $G:(0,\infty)\to(0,\infty)$
such that
\begin{itemize}
\item[(G1)] $G(t)=t$ near $t=\lambda$,
\item[(G2)] $G(t)=1$ for $t\leq\lambda e^{-C}$, and
\item[(G3)] $t(\log G)'(t)\le 1$ for all $t$.
\end{itemize}
Since $\lambda e^{-C}<1$,
one can easily give such a $G$.
Now, set
$H=G\circ K$.
Then, we have the following, as desired.

\begin{Lem}\label{H}
$H$ satisfies the conditions 
$(${\rm H1}$)$ to $(${\rm H4}$)$.
\end{Lem}

\noindent\textbf{Proof.}
Clearly $H$ satisfies (H1).
The validity of (H4) for $H$ follows from that for $K$ 
because by (G1) we have $H=K$ in a neighborhood of 
$\mathcal A\times\{0\}$.
(H2) follows from Lemma \ref{KeC} and (G2).
As stated in the proof of Lemma \ref{KKK}, in order to show (H3)
it suffices to prove that
$\sum_j r_j(\log H)_{r_j}<2$ outside $\mathcal A\times\{0\}$.
We have
$$\sum_j r_j(\log H)_{r_j}
=\sum_j r_j(\log G(K))_{r_j} 
=\sum_j r_j(\log K)_{r_j}K(\log G)'(K).
$$
This is indeed less than $2$,
because
we have already seen in the proof of Lemma \ref{KKK} 
that $\sum_j r_j(\log K)_{r_j}<2$
and  by (G3) we have that $K(\log G)' (K)\le1$.
This verifies (H3).
\hfill $\Box$

\section{Realizing invariant sets in Reeb flows}

\noindent\textbf{Proof of Theorem \ref{main}.}
Suppose that $V$ is a vector field on $T^n$
satisfying  $\sum_{j=1}^nd\theta_j(V)>0$,
and that $\mathcal A\subset T^n$ is a compact invariant set of
the flow generated by $V$.
Then, if we write $V=\sum_j\nu_j\partial_{\theta_j}$,
the sum $\sum_\ell\nu_\ell$ of the coefficient functions is positive.
Take any constants $C$ and $\lambda$ such that
 $0<C<7/9$ and $1<\lambda<e^C$,
and define $C^\infty$ functions $k_j$ and $f$ on $T^n$
by
$k_j=\dfrac{\nu_j}{\sum_\ell\nu_\ell}$
and
$f=\dfrac{2\lambda}{\sum_\ell\nu_\ell}$.
Then, we have
$\sum_j k_j=1$ and 
$fV=2\lambda\sum_j k_j\partial_{\theta_j}$.
By using these $\lambda$ and $k_j$,
we can construct $H$ satisfying the conditions (H1) to (H4)
just as in \S 4. 
If we write $X$ for the contact vector field
corresponding with $H$
with respect to $\alpha_{\mathrm{st}}$
and put $\alpha=\alpha_{\mathrm{st}}/H$,
we see from Lemma \ref{X and H} that 
$X$ is the Reeb vector field of $\alpha$
and satisfies (X1) to (X4).
All the conclusions of Theorem \ref{main}
now follow immediately.
\hfill $\Box$
\vspace{5mm}

\noindent\textbf{Proof of Corollary \ref{Denjoy}.}
In \cite{H}, Hall constructed a $C^\infty$ embedding, say $g$,
of a compact annulus $A$ into itself admitting a Cantor set
as minimal invariant set (Cantor minimal set, for short). 
If we choose an embedding $\iota:A\to T^2$
and extend $\iota\circ g\circ\iota^{-1}:\iota(A)\to\iota(A)$ appropriately,
we obtain a $C^\infty$ diffeomorphism $h_2$ of $T^2$
which is isotopic to the identity and admits a Cantor minimal set.
For any $n\ge3$,
a diffeomorphism $h_{n-1}$ of $T^{n-1}=T^2\times T^{n-3}$
having a Cantor minimal set
is defined by $h_{n-1}=h_2\times\mathrm{id}$.
Now, we consider the suspension flow of $h_{n-1}$.
It is a flow on $T^n$ which has a Denjoy minimal set
and whose orbits are transverse to the fibers of the natural projection
$T^n\to T^{n-1}$.
Passing to a conjugate flow $\psi$ via a suitable linear automorphism of $T^n$
we may assume that the vector field $V$ associated to $\psi$ are transverse to
the plane field $\sum_{j=1}^nd\theta_j=0$.
Thus, we can apply Theorem \ref{main} to $V$ and obtain the conclusion.
\hfill $\Box$

\end{document}